\documentclass[11pt,a4paper]{amsart}
\usepackage[margin=27mm]{geometry}
\usepackage{amsmath,amssymb,amsthm,mathtools}
\usepackage{enumitem,xcolor,fancyhdr,microtype}
\usepackage[colorlinks=true,linkcolor=blue!45!black,citecolor=blue!45!black,urlcolor=blue!50!black]{hyperref}
\definecolor{ink}{RGB}{31,49,65}

\usepackage{tikz}
\usetikzlibrary{arrows.meta}

\newtheorem{theorem}{Theorem}[section]
\newtheorem{proposition}[theorem]{Proposition}
\newtheorem{lemma}[theorem]{Lemma}
\theoremstyle{definition}
\newtheorem{assumption}[theorem]{Assumption}
\theoremstyle{remark}

\numberwithin{equation}{section}
\setlist{topsep=4pt,itemsep=3pt,leftmargin=1.5em}

\newcommand{\QOT}{\texttt{QOT}}
\newcommand{\R}{\mathbb R}
\newcommand{\eps}{\varepsilon}
\newcommand{\dd}{\,\mathrm d}
\newcommand{\ind}{\mathbf 1}

\newcommand{\norm}[1]{\left\lVert#1\right\rVert}

\newcommand{\half}{\tfrac12}

\DeclareMathOperator{\spt}{spt}
\DeclareMathOperator{\interior}{int}

\DeclareMathOperator{\dist}{dist}
\DeclareMathOperator{\Lip}{Lip}

\begin{document}
\title{Support monotonicity for quadratically regularized optimal transport}
\author{Johannes Wiesel}
\address{Johannes Wiesel, Department of Mathematics,
University of Copenhagen}
\email{wiesel@math.ku.dk}

\author{Xingyu Xu}
\address{Xingyu Xu, Department of Electrical and Computer Engineering,
Carnegie Mellon University}
\email{xingyuxu@andrew.cmu.edu}

\begin{abstract}
We study the support of the optimizer $\pi_\eps$ of the
quadratically regularized optimal transport problem for
quadratic cost and compactly supported marginals.
As $\eps$ tends to zero, this support converges to the
support of the Monge coupling. We investigate whether
the regularized supports are nested as $\eps$ decreases.
Under regularity assumptions on the marginal densities,
their domains, and the Brenier map, we establish interior
support inclusions in every dimension. More precisely,
for each compact $K\subset\operatorname{int}X$ and fixed
$\theta\in(0,1)$, we prove that
$\spt\pi_\delta\cap(K\times Y)
\subseteq\spt\pi_\eps\cap(K\times Y)$
whenever $0<\delta\le\theta\eps$ and $\eps$ is sufficiently
small. On these points, the rescaled density $\eps d\pi_\eps/d(\mu\otimes\nu)$ is bounded
below by a positive multiple of $\eps^{2/(d+2)}$.
Our analysis combines sharp potential estimates with
affine rescaling and a Liouville theorem that identifies
the unique limiting density of the rescaled problem.
\end{abstract}

\maketitle

\section{Introduction and main result}

Let $X,Y\subset \R^d$. We study the \emph{quadratically regularized optimal transport} (QOT) problem with quadratic cost: for $\eps>0$, we consider
\begin{equation}\label{eq:primal}
 \QOT_\eps(\mu,\nu)
 =\min_{\pi\in\Pi(\mu,\nu)}
 \left\{\int\frac{|x-y|^2}{2}\dd\pi
 +\frac\eps2\int\left(\frac{\dd\pi}{\dd(\mu\otimes\nu)}\right)^2
                   \dd(\mu\otimes\nu)\right\},
\end{equation}
where $\mu,\nu$ are probability measures on $\R^d$ and $\Pi(\mu,\nu)$ is the set of couplings of $\mu$ and $\nu$.
If $\QOT$ is finite, then it admits a unique optimizer, which we denote by $\pi_\eps$.  Furthermore, the density of $\pi_\eps$ has the following form:
\begin{equation}\label{eq:density}
 \frac{\dd\pi_\eps}{\dd(\mu\otimes\nu)}
       =\eps^{-1}(H_\eps)_+, \qquad\text{where } H_\eps(x,y)=f_\eps(x)+g_\eps(y)-\half|x-y|^2
\end{equation}
for the continuous dual potentials $f_\eps,g_\eps:\R^d\to \R$. We refer to \cite{Nutz2025} for an introduction to $\QOT$ and for a proof of these elementary facts. Throughout this note, we call $H_\eps$ the \emph{slack}. The sign of $H_\eps$ will decide whether a pair belongs to the support. In particular we have
\begin{equation*}
\spt\pi_\eps=\{(x,y)\in X\times Y: H_\eps(x,y)\ge0\},
\end{equation*}
see Lemma \ref{lem:support}.

\begin{figure}[h!]
\begin{center}
\begin{tikzpicture}[
  x=1cm,y=1cm,
  font=\small,
  >=Stealth,
  line cap=round,
  line join=round,
  declare function={
    qotx(\t)=0.6+6.2*\t;
    transport(\t)=0.65+0.5*\t+4*(3*\t*\t-2*\t*\t*\t);
    widthshape(\t)=0.4+0.6*sin(180*\t);
  }
]
  \definecolor{qotouter}{HTML}{E1EEFA}
  \definecolor{qotmiddle}{HTML}{A7CBE9}
  \definecolor{qotinner}{HTML}{5A9DD0}
  \definecolor{qotouteredge}{HTML}{7CAAD0}
  \definecolor{qotmiddleedge}{HTML}{508CB9}
  \definecolor{qotinneredge}{HTML}{2B6E9F}

  \foreach \radius/\shade/\edge in {%
    0.66/qotouter/qotouteredge,
    0.41/qotmiddle/qotmiddleedge,
    0.19/qotinner/qotinneredge%
  }{
    \path[fill=\shade,draw=\edge,line width=0.45pt]
      plot[domain=0:1,samples=101,variable=\t]
        ({qotx(\t)},
         {transport(\t)+\radius*widthshape(\t)})
      -- plot[domain=1:0,samples=101,variable=\t]
        ({qotx(\t)},
         {transport(\t)-\radius*widthshape(\t)})
      -- cycle;
  }

  \draw[black,line width=1.05pt]
    plot[domain=0:1,samples=101,variable=\t]
      ({qotx(\t)},{transport(\t)});
  \node[anchor=west] at (7.02,5.15) {$y=T(x)$};

  \draw[->,line width=0.55pt] (0,0) -- (7.75,0)
    node[anchor=west] {$x$};
  \draw[->,line width=0.55pt] (0,0) -- (0,5.95)
    node[anchor=south] {$y$};
  \draw[line width=0.45pt] (0.6,0.055) -- (0.6,-0.055)
    node[below=3pt] {$a$};
  \draw[line width=0.45pt] (6.8,0.055) -- (6.8,-0.055)
    node[below=3pt] {$b$};

  \foreach \height/\shade/\edge/\index in {%
    5.25/qotouter/qotouteredge/1,
    4.77/qotmiddle/qotmiddleedge/2,
    4.29/qotinner/qotinneredge/3%
  }{
    \path[fill=\shade,draw=\edge,line width=0.45pt]
      (0.78,{\height-0.085}) rectangle (1.25,{\height+0.085});
    \node[anchor=west] at (1.40,\height)
      {$\operatorname{spt}\pi_{\varepsilon_\index}$};
  }
  \node at (4.95,5.88)
    {$\varepsilon_1>\varepsilon_2>\varepsilon_3>0$};

\end{tikzpicture}
\end{center}
\caption{Illustration of nested supports for $d=1,$ where $\operatorname{spt}\pi^\ast\subseteq\operatorname{spt}\pi_{\varepsilon_3}
\subseteq\operatorname{spt}\pi_{\varepsilon_2}\subseteq\operatorname{spt}\pi_{\varepsilon_1}$}\label{fig:1}
\end{figure}
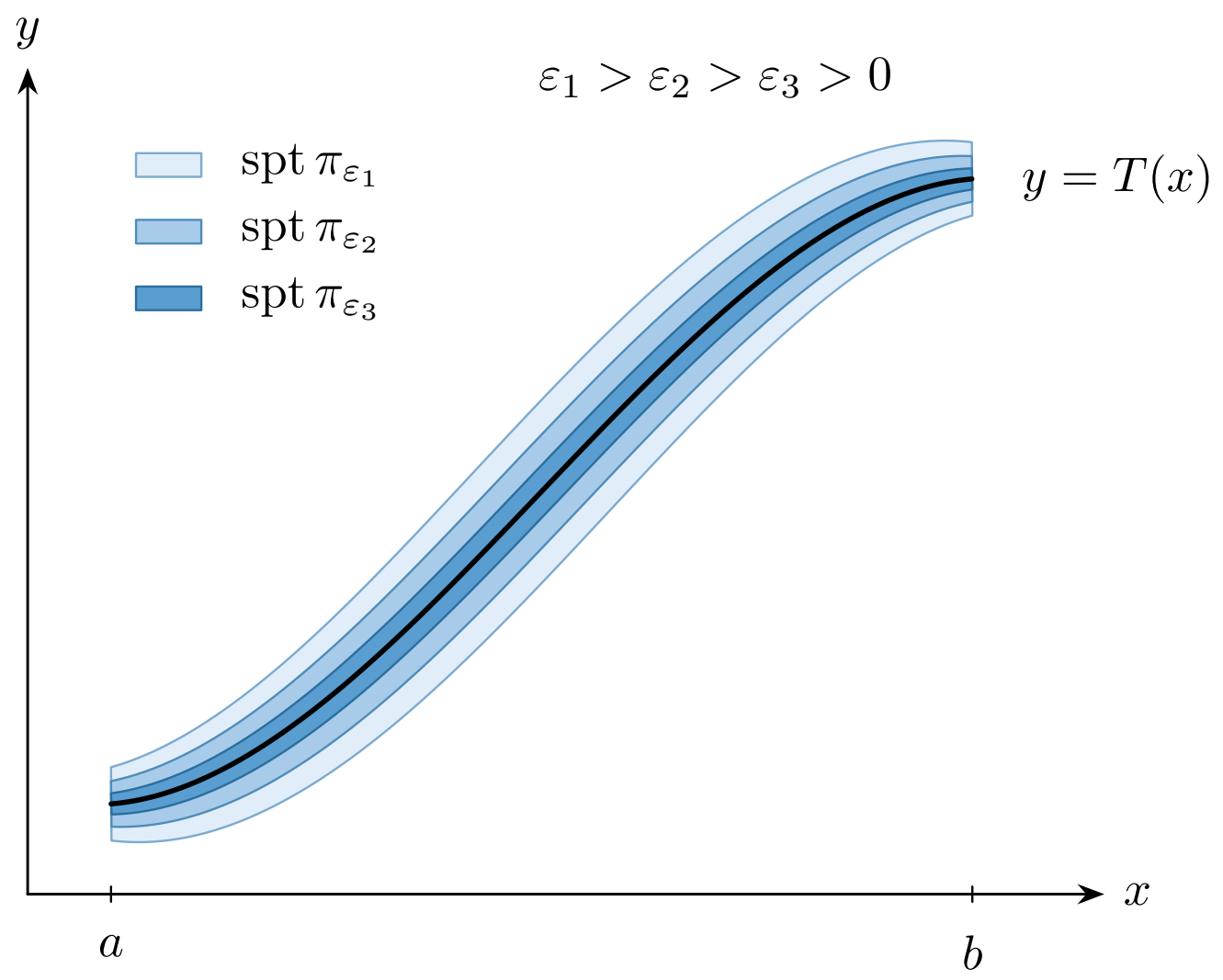

Under regularity assumptions it is well-known that $\spt \pi_\eps$ converges in Hausdorff distance to the support of the OT-optimal coupling $\pi^\ast\in \Pi(\mu,\nu)$; see \cite[Theorem 4.1]{Nutz2025}. This motivates the two complementary questions: 
\begin{enumerate}
\item What is the optimal rate of this convergence?
\item Are the supports eventually nested? More precisely, does there exist $\eps_0>0$ such that $\spt\pi_\delta \subseteq \spt\pi_\eps$ for $\delta\le \eps\le \eps_0$?
\end{enumerate}
Taken together, these questions ask about the main characteristics of the convergence of $\spt \pi_\eps$ to $\spt \pi^\ast$: quantitative rates provide localization bounds for the regularized supports around $\spt\pi^\ast$. An affirmative
answer to the second question ensures that these supports form nested outer approximations of $\spt\pi^\ast$: each contains the limiting support, while their excess shrinks as the regularization decreases (see Figure \ref{fig:1}). 

Question~(1) has been addressed in \cite{gonzalez2026sparsity} in one dimension and in \cite{wiesel2025sparsity} in general dimensions, with sharp rates for the self-transport case. Sharp interior sparsity bounds are obtained in \cite{GonzalezSanzGvalaniKoch2026LocalSparsity}. More recently, \cite{GonzalezSanzNutz2026GeometryI}
establish sharp global convergence estimates, while \cite{GonzalezSanzNutz2026GeometryII} derive detailed interior and boundary asymptotics. 

On the other hand, relatively little progress has been made on Question~(2): A systematic analysis of support
monotonicity in discrete QOT is developed in \cite{gonzalez2024monotonicity}, which provides a geometric characterization of monotonicity, positive results, and explicit counterexamples to nesting along the full regularization path. For fixed finitely supported marginals, however, the regularized optimizer is constant for sufficiently small regularization; see \cite{GonzalezSanzNutz2025LinearPrograms}. Thus these counterexamples do not rule out eventual nesting in the small-regularization regime. To the best of our knowledge, an answer to Question~(2) has not been established for continuous marginals. Our paper thus constitutes the first attempt to fill this gap: we establish support-inclusion results in this regime for quadratic cost, under
regularity assumptions on the marginal densities and the Brenier map.

Before stating our main result, we specify the standing assumptions used throughout the paper. These coincide with \cite[Assumption~1.1]{GonzalezSanzNutz2026GeometryI}. For clarity, we also state explicitly the resulting regularity and ellipticity properties of the Brenier potentials; see \cite[Lemma~2.1]{GonzalezSanzNutz2026GeometryI}.

\begin{assumption}\label{ass:data}
Let $d\ge1$ and $0<\alpha<1$. The sets $X,Y$ are closures of bounded uniformly convex domains with $C^{3,\alpha}$ boundaries; for $d=1$, they are nondegenerate compact intervals. The probability densities are $C^{1,\alpha}(X)$ and $C^{1,\alpha}(Y)$ respectively, and satisfy
\begin{align}\label{eq:density_bounds}
\mu=\rho_\mu\ind_X\dd x,\qquad \nu=\rho_\nu\ind_Y\dd y,
\qquad 0<m\le\rho_\mu,\rho_\nu\le M<\infty.
\end{align}
The Brenier potential $\varphi$ and its conjugate $\varphi^*$ are $C^3$ up to their respective boundaries. The maps $T=\nabla\varphi:X\to Y$ and $S=\nabla\varphi^*:Y\to X$ are inverse diffeomorphisms, and
\begin{align}\label{eq:hessian}
a_0I\preceq \nabla^2\varphi,\nabla^2\varphi^*\preceq b_0I,
\qquad
\norm{\nabla^3\varphi}_\infty+\norm{\nabla^3\varphi^*}_\infty\le M_3,
\end{align}
i.e.~all eigenvalues are in between $a_0$ and $b_0$, where $0<a_0\le b_0<\infty$. Derivatives at the boundary are continuous extensions of interior derivatives.
\end{assumption}

To simplify notation we set
\begin{equation}\label{eq:constants}
\ell_\eps := \eps^{1/(d+2)},\qquad \gamma := \frac2{d+2},\qquad
c_d := \frac12\left(\frac{d+2}{\omega_d}\right)^\gamma,
\end{equation}
where $\omega_d$ is the volume of the unit ball in $\R^d$. Our main result is the following:

\begin{theorem}\label{thm:main}
Under Assumption~\ref{ass:data}, let $K\subset\interior X$ be compact and fix $0<\theta<1$. There exists $\eps_{K,\theta}>0$ such that for all
\[
0<\eps<\eps_{K,\theta},\quad 0<\delta\le\theta\eps,
\quad (x,y)\in\spt\pi_\delta,\quad x\in K
\]
we have
\begin{equation}\label{eq:main}
H_\eps(x,y)\ge
\frac{c_d}{2}M^{-\gamma}(1-\theta^\gamma)\eps^\gamma>0.
\end{equation}
Consequently,
\[
\spt\pi_\delta\cap(K\times Y)
\subseteq\spt\pi_\eps\cap(K\times Y).
\]
\end{theorem}

In particular, the threshold $\eps_{K,\theta}>0$ depends on $K$, $\theta$, whereas the prefactor in the lower bound is independent of $K$. As $K$ is compact, we have $\dist(K,\partial X)>0$. Thus, in the rescaled coordinates used in the proof, the boundary recedes to infinity uniformly over $K$. 

We expect that Assumption \ref{ass:data} can be weakened. The current setting allows us to directly use results from \cite{GonzalezSanzNutz2026GeometryI}, especially the sharp potential bounds in their Theorem~1.3. We also expect a global version of the support inclusion:
for all $\delta\le \theta\eps$, $\spt\pi_\delta\subseteq\spt\pi_\eps$ and $\eps$ is sufficiently small, with no restriction
to $K\times Y$. However, we currently only have a proof for this general result in dimensions $d>200$, as the additional difficulty is to control points whose distance from the boundary is comparable to the shrinking length $\ell_\eps$. Here we present the interior argument, which applies in every dimension. Establishing the global extension in the remaining dimensions, possibly through a rigorous
computer-assisted analysis, is left for future work. Similarly we conjecture that there exists $\eps_0>0$ such that
the support inclusion holds whenever
$0<\delta\le\eps\le\eps_0$.

We now provide a quick roadmap of the proof of Theorem \ref{thm:main}. For this we use the notation
\begin{align}\label{eq:notation}
f_*(x):=\frac12|x|^2-\varphi(x),\quad
 g_*(y):=\frac12|y|^2-\varphi^*(y),\quad
r_\eps:=f_\eps-f_*,\quad s_\eps:=g_\eps-g_*.
\end{align}
We note that 
\begin{equation}\label{eq:decomposition}
H_\eps=r_\eps+s_\eps-D,\qquad \text{where }
D(x,y):=\varphi(x)+\varphi^*(y)-\langle x, y \rangle
\end{equation}
is independent of $\eps$. Thus, to compare $H_\eps$ and $H_\delta$, it suffices to estimate the correction sum $r_\eps(x)+s_\eps(y)$ accurately at both scales. 

We bound the distance of the regularized support from the graph of the Brenier map and obtain Lipschitz bounds for the differences between the regularized and unregularized potentials in Section~\ref{sec:prelim}. In Section~\ref{sec:rescaling} we then introduce affine coordinates and rescale these differences. Compactness allows us to pass to the limit in the marginal equations, obtaining equations on $\R^d$. In Section~\ref{sec:classification} we establish an exact comparison and a Liouville theorem for bounded solutions, that identify the unique limiting slack. Compactness and uniqueness then yield convergence uniformly over $K$, see Section~\ref{sec:uniform}. Finally, we compare the slack at the
two regularization parameters $\delta$ and $\eps$ in Section~\ref{sec:comparison}. Appendix \ref{app:liouville} contains a detailed proof of the Liouville theorem.

\textbf{Related Literature.}
Early contributions to the $\QOT$ problem include
\cite{BlondelSeguyRolet2018} in a discrete setting
and \cite{EssidSolomon2018} for transport on graphs.
For continuous marginals, \cite{LorenzMannsMeyer2021}
establish existence and strong duality and develop numerical
methods. A broader framework covering quadratic regularization
together with convergence results under smoothing and
discretization, is developed in \cite{LorenzMahler2022}.
Quantitative convergence rates for regularized transport costs
are obtained via quantization in
\cite{EcksteinNutz2024Quantization}, while
\cite{GarrizMolinaGonzalezSanzMordant2024Infinitesimal}
derive exact first-order asymptotics for the QOT cost,
revealing a connection with the porous medium equation.
For related regularity results, we also refer to
\cite{GvalaniKoch2026UniformRegularity, ChengKochXiao2026GlobalRegularity}.

\section{Preliminary results}\label{sec:prelim}

Let us start by recalling a few additional results from \cite{Nutz2025}: $f_\eps, g_\eps$ are only unique up to additive shifts, that do not change $H_\eps$. Furthermore, the marginal constraints imply the pointwise constraints
\begin{equation}\label{eq:marginals}
\int_X [H_\eps(x,y)]_+\dd \mu(x)=\eps,
\qquad
\int_Y [H_\eps(x,y)]_+\dd\nu( y)=\eps.
\end{equation}

\begin{lemma}[Concavity of sections]\label{lem:support}
For each fixed variable, $H_\eps$ is concave in the other variable. In particular we have
\begin{equation}\label{eq:support}
\spt\pi_\eps=\overline{\{H_\eps>0\}}=\{H_\eps\ge0\}.
\end{equation}
in the relative topology on $X\times Y$.
\end{lemma}

\begin{proof}
Set $u_\eps(x):=|x|^2/2-f_\eps(x)$ and $v_\eps(y):=|y|^2/2-g_\eps(y)$. Then $H_\eps(x,y)=\langle x, y \rangle-u_\eps(x)-v_\eps(y)$. The function
\[
\Psi(x,a)=\int_Y[\langle x, y \rangle-a-v_\eps(y)]_+\dd\nu(y)
\]
is convex in $(x,a)$, because its integrand is the positive part of an affine function of those variables. It is strictly decreasing in $a$ whenever its value is positive. Indeed, if $a_1<a_2$ and $\Psi(x,a_2)>0$, decreasing $a$ raises the integrand by $a_2-a_1$ on a set of positive measure. Since $\Psi(x,u_\eps(x))=\eps>0$ by \eqref{eq:marginals}, its sublevel set $\{\Psi\le\eps\}$ is exactly the epigraph $\{a\ge u_\eps(x)\}$. That epigraph is convex, hence $u_\eps$ is convex. A symmetric argument yields convexity of $v_\eps$. Thus $H_\eps =\langle x, y \rangle-u_\eps(x)-v_\eps(y)$ is separately concave.

For \eqref{eq:support}, by virtue of \eqref{eq:density} we have $\{H_\eps > 0\} \subseteq \spt \pi_\eps \subseteq \{H_\eps \ge 0\}$; thus it suffices to show $\{H_\eps\ge0\}=\overline{\{H_\eps>0\}}$.
Assume $H_\eps(x, y) = 0$ for some $x \in X$, $y \in Y$. The identity \eqref{eq:marginals} guarantees the existence of some $y' \in Y$ such that $H_\eps(x, y') > 0$. Concavity ensures $H_\eps$ is positive in the interior of the segment connecting $(x, y)$ and $(x, y')$. This implies $(x, y) \in \overline{\{H_\eps>0\}}$ by continuity, hence $\{H_\eps\ge0\}=\overline{\{H_\eps>0\}}$ as desired.
\end{proof}

\begin{lemma}[Sharp global potential convergence, cf.~{\cite[Theorem 1.3]{GonzalezSanzNutz2026GeometryI}}]\label{thm:external-height}
Under Assumption~\ref{ass:data}, fix $x_0\in X$ and choose the additive constant of $f_\eps$ so that $r_\eps(x_0)=0$. There are $C_0,\eps_0>0$, independent of $\eps$, such that
\begin{equation}\label{eq:height}
\norm{r_\eps}_\infty+\norm{s_\eps}_\infty\le C_0\ell_\eps^2
\qquad(0<\eps<\eps_0).
\end{equation}
\end{lemma}

The normalization of the potentials in the above lemma is used throughout the remaining paper.
Recall $D(x,y)=\varphi(x)+\varphi^*(y)-\langle x, y \rangle$.
\begin{lemma}[Quadratic bounds for the transport gap]\label{lem:gap}
For all $(x,y)\in X\times Y$,
\[
\frac{a_0}{2}|y-T(x)|^2\le D(x,y)\le\frac{b_0}{2}|y-T(x)|^2,
\qquad
D(x,y)\ge\frac{a_0}{2}|x-S(y)|^2.
\]
\end{lemma}
\begin{proof}
Put $q=T(x)$ and $z=y-q$. Conjugacy gives $\varphi(x)+\varphi^*(q)=\langle x,q\rangle $ and $\nabla\varphi^*(q)=x$. Consequently,
\[
D(x,y)=\varphi^*(y)-\varphi^*(q)-\langle \nabla\varphi^*(q),y-q\rangle
=\int_0^1(1-t)z^\top \nabla^2\varphi^*(q+tz)z\dd t.
\]
The segment $\{q+tz: t\in[0,1]\}$ lies in the convex set $Y$. Using \eqref{eq:hessian} and $\int_0^1(1-t)\dd t=1/2$ yields the first claim. Interchanging $x,y$ and $\varphi,\varphi^*$ proves the second claim.
\end{proof}

\begin{proposition}[Uniform localization]\label{prop:localization}
There are constants $C_{\mathrm{loc}},C_1$ such that, for small $\eps$,
\begin{align}
(x,y)\in\spt\pi_\eps&\ \Rightarrow\
|y-T(x)|+|x-S(y)|\le C_{\mathrm{loc}}\ell_\eps,\label{eq:localization}\\
\Lip(r_\eps)+\Lip(s_\eps)&\le C_1\ell_\eps.\label{eq:slopes}
\end{align}
Here $\Lip(a)$ is the least $L$ such that $|a(z)-a(z')|\le L|z-z'|$ on its domain.
\end{proposition}
\begin{proof}
We first prove \eqref{eq:localization}: at a support point we have $H_\eps\ge0$, so \eqref{eq:decomposition} and \eqref{eq:height} give $D\le r_\eps+s_\eps\le C_0\ell_\eps^2$. By Lemma~\ref{lem:gap} we have $D\ge a_0|y-T(x)|^2/2$, so that $|y-T(x)|\le \sqrt{2C_0/a_0}\,\ell_\eps$. The same bound holds for $|x-S(y)|$. This proves \eqref{eq:localization}.

We turn to the proof of \eqref{eq:slopes}. For this we fix $x, x'\in X$ and let $b=x'-x$, $e=|b|$. For $y\in Y$ such that $(x,y)\in \spt\pi_\eps$ we have
\[
D(x',y)-D(x,y) = \varphi(x') -\varphi(x)-\langle x'-x,y\rangle
=\int_0^1 \langle T(x+tb)-y, b\rangle \dd t.
\]
As $T(x+tb)-y = [T(x+tb)-T(x)] + [T(x)-y]$, \eqref{eq:hessian} and \eqref{eq:localization} yield
\begin{align}\label{eq:help}
|D(x',y)-D(x,y)|
\le C_{\mathrm{loc}}\ell_\eps e+b_0e^2.
\end{align}
We now show 
\begin{align}\label{eq:toshow}
r_\eps(x')-r_\eps(x)\le C_{\mathrm{loc}}\ell_\eps e+b_0e^2
\end{align}
In fact if $r_\eps(x')-r_\eps(x)\ge C_{\mathrm{loc}}\ell_\eps e+b_0e^2+a$ for some $a>0$, then 
\[
H_\eps(x',y)-H_\eps(x,y)= r_\eps(x')-r_\eps(x) -[D(x',y)-D(x,y)] \stackrel{\eqref{eq:help}}{\ge} a
\]
for every $y\in E:=\{y\in Y: H_\eps(x,y)\ge 0\}$. Then 
$$\eps\stackrel{\eqref{eq:marginals}}{=}\int_Y [H_\eps(x',y)]_+\nu(\dd y)\ge \int_E ([H_\eps(x,y)]_+ +a)\nu(\dd y) \ge \eps+a\nu(E)>\eps$$ as $\nu(E)>0$ by \eqref{eq:marginals}, a contradiction. Thus \eqref{eq:toshow} follows.
Interchanging $x$ and $x'$ gives the same bound for the negative difference. Now we subdivide the segment from $x$ to $x'$ into $N$ equal pieces and sum these estimates. Convexity of $X$ ensures that all subdivision points lie in $X$, and thus
\[
|r_\eps(x')-r_\eps(x)|
\le C_{\mathrm{loc}}\ell_\eps e+\frac{b_0e^2}{N}.
\]
Letting $N\to\infty$ proves $\Lip(r_\eps)\le C_{\mathrm{loc}}\ell_\eps$. The estimate for $s_\eps$ follows in the same way.
\end{proof}

\section{Local coordinates and full-space limit}\label{sec:rescaling}

In this section we use the affine normalization described in
\cite[Section~2.4]{GvalaniKoch2026UniformRegularity};
see also
\cite[Section~3.3]{ChengKochXiao2026GlobalRegularity}.
We choose the coordinates to normalize the derivative
of the Brenier map at the base point. 

In more detail, we now fix a base point $p\in K$ and set $q:=T(p)$ as well as
\begin{equation}\label{eq:affine-data}
J:=\nabla T(p),\qquad \lambda(p):=\sqrt{\rho_\mu(p)\rho_\nu(q)},\qquad
h=h_\eps(p):=\left(\frac{\eps}{\lambda(p)}\right)^{1/(d+2)},
\end{equation}
where we recall that $\rho_\mu, \rho_\nu$ are the densities of $\mu$ and $\nu$ respectively.
The matrix $J$ is symmetric positive definite. Throughout, matrix norms are operator norms, i.e.~$\norm{L}=\sup_{|z|=1}|Lz|$ for a matrix $L$. Recalling \eqref{eq:density_bounds} and \eqref{eq:hessian},
\begin{align}\label{eq:norm_bound}
\norm{J^{1/2}}\le\sqrt{b_0},\qquad
\norm{J^{-1/2}}\le a_0^{-1/2},\qquad
m\le\lambda(p)\le M.
\end{align}
The numbers $h$ and $\ell_\eps$ defined in \eqref{eq:constants} satisfy
\begin{align}\label{eq:comparison:hl}
h=\ell_\eps\lambda(p)^{-1/(d+2)}.
\end{align}

We now introduce the rescaled coordinates $\xi, \eta\in \R^d$ via 
\begin{equation}\label{eq:coordinates}
x=p+hJ^{-1/2}\xi,\qquad y=q+hJ^{1/2}\eta.
\end{equation}
These describe displacements from $p,q$ in units of $h$, after adjusting the quadratic forms by $J$. In particular, the base pair $(p,q)$ has coordinates $(0,0)$.

\begin{lemma}\label{lem:rescaling}
Fix $L>0$ and take $\varepsilon$ sufficiently small. Then for all $|\xi|,|\eta|\le L$, uniformly over $p\in K$, we have
\begin{equation}\label{eq:gap-taylor}
h^{-2}D(x,y)=h^{-2}D(p+hJ^{-1/2}\xi,q+hJ^{1/2}\eta)
=\frac12|\xi-\eta|^2+O_L(h).
\end{equation}
The transformed measures are
\begin{align}
\dd\mu(x)&=\lambda(p)h^d
\frac{\rho_\mu(p+hJ^{-1/2}\xi)}{\rho_\mu(p)}\dd\xi,\label{eq:source-density}\\
\dd\nu(y)&=\lambda(p)h^d
\frac{\rho_\nu(q+hJ^{1/2}\eta)}{\rho_\nu(q)}\dd\eta.\label{eq:target-density}
\end{align}
We have $$\lim_{\eps\downarrow 0}\sup_{|\xi|\le L, p\in K} \Big|\frac{\rho_\mu(p+hJ^{-1/2}\xi)}{\rho_\mu(p)}-1\Big|=0=\lim_{\eps\downarrow 0}\sup_{|\eta|\le L, p\in K} \Big|\frac{\rho_\nu(q+hJ^{1/2}\eta)}{\rho_\nu(q)}-1\Big|.$$
\end{lemma}

\begin{proof}
By definition we have $D(p,q)=0$, $\nabla_xD(p,q)=T(p)-q=0$, and $\nabla_yD(p,q)=S(q)-p=0$. Differentiating $S(T(p))=p$ gives $\nabla S(q)J=I$, so $\nabla^2\varphi^*(q)=J^{-1}$. Taylor expansion gives, with $a=x-p$ and $b=y-q$,
\[
D(x,y)=\frac12a^\top Ja+\frac12b^\top J^{-1}b-\langle a, b\rangle 
+O(|a|^3+|b|^3).
\]
By \eqref{eq:coordinates},
\[
a^\top Ja=h^2|\xi|^2,\quad
b^\top J^{-1}b=h^2|\eta|^2,\quad
\langle a, b\rangle =h^2\xi^\top J^{-1/2}J^{1/2}\eta=h^2\langle \xi,\eta\rangle .
\]
Thus \eqref{eq:gap-taylor} follows. Recalling \eqref{eq:hessian} and \eqref{eq:norm_bound}, the remainder in the Taylor expansion of $D$ is bounded by $C_Lh^3$, uniformly over $p\in K$. Dividing by $h^2$ thus proves \eqref{eq:gap-taylor}.

The change-of-variables formula for $T_\#\mu=\nu$ gives
\[
\rho_\mu(p)=\rho_\nu(q)\det J.
\]
Since $\det J>0$,
\[
\rho_\mu(p)(\det J)^{-1/2}
=\rho_\nu(q)(\det J)^{1/2}
=\sqrt{\rho_\mu(p)\rho_\nu(q)}=\lambda(p).
\]
The Jacobians in  \eqref{eq:coordinates} are $h^d(\det J)^{-1/2}$ and $h^d(\det J)^{1/2}$, proving \eqref{eq:source-density}--\eqref{eq:target-density}. If $|\xi|,|\eta|\le L$, then $|x-p|,|y-q|$ are $O_L(h)$, as $$\sup_{p\in K}h_\eps(p)\le \Big( \frac{\eps}{m}\Big)^{1/(d+2)}\to 0.$$ Uniform continuity of the densities, and their positive lower bound, make the ratios tend uniformly to one, even when $p$ varies.
\end{proof}

Recalling \eqref{eq:decomposition} we now characterise the limits of the rescaled integral equations \eqref{eq:marginals} with the help of Lemma \ref{lem:rescaling}.

\begin{proposition}\label{prop:compactness}
Let $\eps_j\downarrow0$ and $p_j\in K$. Define $q_j,J_j,\lambda_j,h_j$ for $\eps=\eps_j$ as in \eqref{eq:affine-data}. Define the rescaled corrections by
\begin{align}
A_j(\xi)&:=h_j^{-2}r_{\eps_j}(p_j+h_jJ_j^{-1/2}\xi),\\
B_j(\eta)&:=h_j^{-2}s_{\eps_j}(q_j+h_jJ_j^{1/2}\eta),
\end{align}
where we recall $r_{\eps_j},s_{\eps_j}$ from \eqref{eq:notation}. These are defined on 
\begin{align*}
X_j&:=h_j^{-1}J_j^{1/2}(X-p_j) \quad \text{and} \quad Y_j:=h_j^{-1}J_j^{-1/2}(Y-q_j),
\end{align*}
respectively.
Then every fixed ball is contained in $X_j\cap Y_j$ for all sufficiently large $j$. After passing to a subsequence, $A_j,B_j$ converge uniformly on every fixed ball to bounded Lipschitz functions $A,B$ on $\R^d$. They satisfy
\begin{align}\label{eq:frozen}
\begin{split}
\int_{\R^d}[A(\xi)+B(\eta)-\tfrac12|\xi-\eta|^2]_+\dd\eta&=1
\quad \text{for every }\xi,\\
\int_{\R^d}[A(\xi)+B(\eta)-\tfrac12|\xi-\eta|^2]_+\dd\xi&=1
\quad \text{for every }\eta.
\end{split}
\end{align}
\end{proposition}

\begin{proof}
We proceed via three steps.

\emph{Step 1: every ball is contained in $X_j\cap Y_j$ for all sufficiently large $j$.}
Recall that $\xi \in X_j$ if and only if $x=p_j+h_j J_j^{-1/2}\xi \in X$. Thus, we have to show that $|\xi|\le L$ implies $x\in X$.
The sets $K$ and $T(K)$ are compact inside $X$ and $Y$, respectively, thus $\text{dist}(K, \partial X), \text{dist}(T(K), \partial Y)>0$. For fixed $L$, $|\xi|\le L$ implies
\[
|x-p_j|
\le h_j a_0^{-1/2}L
<\dist(K,\partial X)
\]
for all sufficiently large $j$, since $h_j\to0$.
Consequently,
\[
x\in B_{\dist(K,\partial X)}(p_j)\subset X.
\]
Thus $B_L(0)\subset X_j$. The corresponding estimate
\[
|y-q_j|\le h_j\sqrt{b_0}L
<\dist(T(K),\partial Y)
\]
proves $B_L(0)\subset Y_j$ for all sufficiently large $j$.

\emph{Step 2: the functions have a convergent subsequence.}
Fix a ball.  Lemma~\ref{thm:external-height} and Proposition~\ref{prop:localization} give
\[
\norm{A_j}_\infty+\norm{B_j}_\infty
\le C\ell_{\eps_j}^2/h_j^2\stackrel{\eqref{eq:comparison:hl}}{\le} C,
\qquad
\Lip(A_j)+\Lip(B_j)\le C\ell_{\eps_j}/h_j\stackrel{\eqref{eq:comparison:hl}}{\le} C.
\]
Thus, Arzel\`a--Ascoli applies on the closure of this fixed ball. We can now choose subsequences on balls of radii $1,2,\ldots$ and then a diagonal subsequence. The limits exist on all of $\R^d$, have the same bounds, and are Lipschitz.

\emph{Step 3: the marginal equations pass to the limit.}
Fix $\xi$ and set
\[
x:=p_j+h_jJ_j^{-1/2}\xi,
\qquad y_j(\eta):=q_j+h_jJ_j^{1/2}\eta
\]
for sufficiently large $j$, such that $x\in X$.
We then compute
\begin{align*}
\eps_j
&\stackrel{\eqref{eq:marginals}}{=}\int_Y [H_{\eps_j}(x,y)]_+\dd\nu(y)\\
&\stackrel{\eqref{eq:decomposition},\eqref{eq:target-density}}{=}\lambda_jh_j^d
\int_{Y_j}
\left[h_j^2A_j(\xi)+h_j^2B_j(\eta)-D(x,y_j(\eta))\right]_+
\frac{\rho_\nu(y_j(\eta))}{\rho_\nu(q_j)}\dd\eta\\
&=\lambda_jh_j^{d+2}
\int_{Y_j}
\left[A_j(\xi)+B_j(\eta)-h_j^{-2}D(x,y_j(\eta))\right]_+
\frac{\rho_\nu(y_j(\eta))}{\rho_\nu(q_j)}\dd\eta.
\end{align*}
Since
$\lambda_jh_j^{d+2}=\eps_j>0$ by \eqref{eq:affine-data}, we obtain
\begin{align}\label{eq:int_goal}
\int_{Y_j}
\left[A_j(\xi)+B_j(\eta)-h_j^{-2}D(x,y_j(\eta))\right]_+
\frac{\rho_\nu(y_j(\eta))}{\rho_\nu(q_j)}\dd\eta=1.
\end{align}
It remains to take limits in \eqref{eq:int_goal}.  For this we first argue that we can restrict to balls. For fixed $\xi$, set
\begin{align*}
F_j(\eta)
&:=\left[A_j(\xi)+B_j(\eta)
-h_j^{-2}D(x,y_j(\eta))\right]_+ \qquad
F(\eta)
&:=\left[A(\xi)+B(\eta)-\tfrac12|\xi-\eta|^2\right]_+,
\end{align*}
and $w_j(\eta):=\rho_\nu(y_j(\eta))/\rho_\nu(q_j)$, where $F_j, w_j$ are defined on $Y_j$.
If $F_j(\eta)>0$, then
$H_{\eps_j}(x,y_j(\eta))>0$ so that $(x,y_j(\eta))\in \spt\pi_{\eps_j}$.
Using \eqref{eq:localization} and $q_j=T(p_j)$ together with \eqref{eq:hessian} we obtain
\begin{align*}
|\eta|
&\le h_j^{-1}a_0^{-1/2}
\bigl(|y_j(\eta)-T(x)|+|T(x)-q_j|\bigr)\\
&\le a_0^{-1/2}
\left(
C_{\mathrm{loc}}\frac{\ell_{\eps_j}}{h_j}
+b_0a_0^{-1/2}|\xi|
\right)\\
&\stackrel{\eqref{eq:comparison:hl}}{\le} C_*(1+|\xi|).
\end{align*}
Set $M_0:=\|A\|_\infty+\|B\|_\infty$. Then
\[
F(\eta)>0
\ \Rightarrow\
\tfrac12|\xi-\eta|^2<A(\xi)+B(\eta)\le M_0
\ \Rightarrow\
|\eta|<|\xi|+\sqrt{2M_0}.
\]
Choose $R>\max\{C_*(1+|\xi|),\,|\xi|+\sqrt{2M_0}\}$.
By the above we have $F_j=0$ on $Y_j\setminus B_R(0)$ and $F=0$ on $\R^d\setminus B_R(0)$.
By Step~1, $B_R(0)\subseteq Y_j$ for all sufficiently large $j$.

Now we can put things together. Indeed, the local uniform convergence of $A_j,B_j$,
\eqref{eq:gap-taylor}, and $|u_+-v_+|\le|u-v|$ imply
\[
\|F_j-F\|_{L^\infty(B_R(0))}\to 0.
\]
Lemma~\ref{lem:rescaling} and \eqref{eq:density_bounds} give
\[
\|w_j-1\|_{L^\infty(B_R(0))}\to 0,
\qquad
\|w_j\|_{L^\infty(B_R(0))}\le M/m.
\]
Consequently, with all norms below taken on $B_R(0)$,
\begin{align*}
\left|\int_{B_R}F_jw_j\dd\eta-\int_{B_R}F\dd\eta\right|
&\le |B_R(0)|\left(
\|w_j\|_\infty\|F_j-F\|_\infty
+\|F\|_\infty\|w_j-1\|_\infty
\right)\to 0.
\end{align*}
Lastly, it follows from \eqref{eq:int_goal} that
\begin{align*}
\int_{\R^d}F(\eta)\dd\eta=\int_{B_R(0)}F(\eta)\dd\eta
&=\lim_{j\to\infty}\int_{B_R(0)}F_j(\eta)w_j(\eta)\dd\eta\\
&=\lim_{j\to\infty}\int_{Y_j}F_j(\eta)w_j(\eta)\dd\eta
=1.
\end{align*}
This proves the first identity in \eqref{eq:frozen}.
The second follows by interchanging the two variables,
using \eqref{eq:source-density} and the estimate
$|x-S(y)|\le C_{\mathrm{loc}}\ell_{\eps_j}$
from \eqref{eq:localization}.
\end{proof}

\section{The full-space solution}\label{sec:classification}

The integral equations \eqref{eq:frozen} allow us to identify all bounded continuous solutions $A,B$ of the limiting marginal equations, up to opposite additive constants. For the proof we need the following Liouville principle, which is proved in Appendix~\ref{app:liouville}.

\begin{lemma}[Liouville principle]\label{lem:liouville}
Let $\mathcal J:R^d\times\R^d\to \R$ be a jointly Borel measurable symmetric and satisfy
\begin{equation}\label{eq:Jbounds}
j_-\ind_{\{|x-z|<r\}}\le\mathcal J(x,z)
\le j_+\ind_{\{|x-z|\le R_0\}},
\qquad 0<r\le R_0,
\end{equation}
with $j_-,j_+>0$. Then every bounded continuous solution of
\begin{align}\label{eq:u}
\int_{\R^d}\mathcal J(x,z)[u(x)-u(z)]\dd z=0
\qquad \forall x\in\R^d
\end{align}
is constant.
\end{lemma}

\begin{proposition}\label{prop:classification}
Every bounded continuous pair $A,B$ satisfying \eqref{eq:frozen} has the form
\[
A(\xi)=c_d/2+c,\qquad B(\eta)=c_d/2-c
\]
for some constant $c$. In particular, its slack is
\begin{equation}\label{eq:cap}
H_0(\xi,\eta):= A(\xi)+B(\eta)-\frac{1}{2}|\xi-\eta|^2=c_d-\frac12|\xi-\eta|^2.
\end{equation}
\end{proposition}

\begin{proof}
Our aim is to apply the Liouville principle in the form of Lemma \ref{lem:liouville} to the difference $\int ([H]_+-[H_0]_+)=0$. For this we have to identify a symmetric kernel. This is what the main effort in this proof is spent on. 
We proceed in four steps.

\emph{Step 1: construct a constant comparison pair.}
For $c>0$, using polar coordinates gives
\begin{align*}
\int_{\R^d}\left(c-\frac12|z|^2\right)_+\dd z
&=d\omega_d\int_0^{\sqrt{2c}}
\left(c-\frac12t^2\right)t^{d-1}\dd t
=\frac{2^{(d+2)/2}\omega_d}{d+2}c^{(d+2)/2}.
\end{align*}
By \eqref{eq:constants}, this expression equals one when $c=c_d$.
Thus, setting
\[
A_0(\xi)=B_0(\eta):=\frac{c_d}{2},
\qquad
H_0(\xi,\eta)=c_d-\frac12|\xi-\eta|^2,
\]
we obtain
\[
\int_{\R^d}[H_0(\xi,\eta)]_+\dd\eta=1
\quad\text{for every }\xi,
\qquad
\int_{\R^d}[H_0(\xi,\eta)]_+\dd\xi=1
\quad\text{for every }\eta.
\]

\emph{Step 2: subtract the marginal identities.}
For $\alpha,\beta\in\R$, the fundamental theorem of calculus gives
\begin{equation}\label{eq:linearization}
[\alpha]_+-[\beta]_+
=(\alpha-\beta)\int_0^1
\ind_{\{(1-t)\beta+t\alpha>0\}}\dd t.
\end{equation}
Set
\begin{align*}
H(\xi,\eta)&:=A(\xi)+B(\eta)-\frac12|\xi-\eta|^2,\\
u(\xi)&:=A(\xi)-\frac{c_d}{2},
\qquad
v(\eta):=B(\eta)-\frac{c_d}{2},
\end{align*}
and define
\[
\mathcal K(\xi,\eta):=
\int_0^1
\ind_{\{(1-t)H_0(\xi,\eta)+tH(\xi,\eta)>0\}}\dd t.
\]
Since $H-H_0=u+v$, \eqref{eq:linearization} yields
\begin{align}\label{eq:H}
[H(\xi,\eta)]_+-[H_0(\xi,\eta)]_+
=\mathcal K(\xi,\eta)\bigl[u(\xi)+v(\eta)\bigr].
\end{align}
We first establish bounds for $\mathcal K$. Set
\[
M_*:=\norm{A}_\infty+\norm{B}_\infty,
\qquad
r_*:=\sqrt{c_d},
\qquad
R_*:=\sqrt{2(M_*+c_d)}.
\]
If $|\xi-\eta|>R_*$, then
\[
H(\xi,\eta)\le M_*-\frac12|\xi-\eta|^2<0,
\qquad
H_0(\xi,\eta)=c_d-\frac12|\xi-\eta|^2<0.
\]
Consequently, $\mathcal K(\xi,\eta)=0$.

If $|\xi-\eta|<r_*$, then
\[
H_0(\xi,\eta)\ge\frac{c_d}{2},
\qquad
H(\xi,\eta)\ge-M_*-\frac{c_d}{2}.
\]
Hence
\[
(1-t)H_0(\xi,\eta)+tH(\xi,\eta)
\ge\frac{c_d}{2}-t(M_*+c_d)>0
\]
whenever
\[
0\le t<k_*:=\frac{c_d}{2(M_*+c_d)}.
\]
We have therefore proved
\begin{align}\label{eq:kbounds}
k_*\ind_{\{|\xi-\eta|<r_*\}}
\le\mathcal K(\xi,\eta)
\le\ind_{\{|\xi-\eta|\le R_*\}}.
\end{align}
Define
\begin{align}\label{eq:H2}
a(\xi):=\int_{\R^d}\mathcal K(\xi,\eta)\dd\eta,
\qquad
b(\eta):=\int_{\R^d}\mathcal K(\xi,\eta)\dd\xi.
\end{align}
The preceding bounds imply
\begin{align}\label{eq:Kmass}
a_-\le a(\xi)\le a_+,
\qquad
a_-\le b(\eta)\le a_+,
\end{align}
where $a_-:=k_*\omega_dr_*^d>0, a_+:=\omega_dR_*^d<\infty.$
Subtracting the marginal identities for $H$ and $H_0$ and recalling \eqref{eq:H} and \eqref{eq:H2} now gives
\begin{align}\label{eq:bipartite}
\begin{split}
0
&=\int_{\R^d}\bigl([H(\xi,\eta)]_+
-[H_0(\xi,\eta)]_+\bigr)\dd\eta=a(\xi)u(\xi)
+\int_{\R^d}\mathcal K(\xi,\eta)v(\eta)\dd\eta,\\
0
&=\int_{\R^d}\bigl([H(\xi,\eta)]_+
-[H_0(\xi,\eta)]_+\bigr)\dd\xi=b(\eta)v(\eta)
+\int_{\R^d}\mathcal K(\xi,\eta)u(\xi)\dd\xi.
\end{split}
\end{align}

\emph{Step 3: derive an equation with a symmetric kernel.}
In Step 2 we saw that $\mathcal K$ is a good candidate for a kernel $\mathcal J$ in Liouville's theorem; however its definition does not guarantee symmetry. This step symmetrises it.

Since $b(\eta)\ge a_->0$, the second identity in
\eqref{eq:bipartite} implies
\[
v(\eta)
=-\frac1{b(\eta)}
\int_{\R^d}\mathcal K(\zeta,\eta)u(\zeta)\dd\zeta.
\]
Define
\[
\mathcal J(\xi,\zeta):=
\int_{\R^d}
\frac{\mathcal K(\xi,\eta)\mathcal K(\zeta,\eta)}
{b(\eta)}\dd\eta.
\]
This kernel is measurable, nonnegative, and symmetric. Moreover, Tonelli's theorem gives
\begin{align}
\int_{\R^d}\mathcal J(\xi,\zeta)\dd\zeta
&=\int_{\R^d}
\frac{\mathcal K(\xi,\eta)}{b(\eta)}
\left[
\int_{\R^d}\mathcal K(\zeta,\eta)\dd\zeta
\right]\dd\eta
=\int_{\R^d}
\frac{\mathcal K(\xi,\eta)}{b(\eta)}
b(\eta)\dd\eta
=a(\xi). \label{eq:tonelli}
\end{align}
In particular,
\[
\int_{\R^d}\mathcal J(\xi,\zeta)|u(\zeta)|\dd\zeta
\le a_+\norm{u}_\infty<\infty.
\]
Thus Fubini's theorem applies when substituting the expression
for $v$ into the first identity in \eqref{eq:bipartite}:
\begin{align*}
a(\xi)u(\xi)
&=-\int_{\R^d}\mathcal K(\xi,\eta)v(\eta)\dd\eta
=\int_{\R^d}\mathcal J(\xi,\zeta)u(\zeta)\dd\zeta.
\end{align*}
Using \eqref{eq:tonelli} 
we conclude that
\[
\int_{\R^d}
\mathcal J(\xi,\zeta)
\bigl[u(\xi)-u(\zeta)\bigr]\dd\zeta=0
\qquad\text{for every }\xi\in\R^d.
\]

\emph{Step 4: apply the Liouville principle.}
Suppose $|\xi-\zeta|<r_*/2$. For every
$\eta\in B_{r_*/2}(\xi)$, we have
\[
|\xi-\eta|<r_*,
\qquad
|\zeta-\eta|
\le|\zeta-\xi|+|\xi-\eta|<r_*.
\]
Hence both factors of $\mathcal K$ are at least $k_*$, and
\eqref{eq:kbounds} together with \eqref{eq:Kmass} gives
\[
\mathcal J(\xi,\zeta)
\ge
\int_{B_{r_*/2}(\xi)}\frac{k_*^2}{a_+}\dd\eta
=\frac{k_*^2\omega_d(r_*/2)^d}{a_+}.
\]
On the other hand, again by \eqref{eq:kbounds}, a nonzero product
$\mathcal K(\xi,\eta)\mathcal K(\zeta,\eta)$ requires
\[
|\xi-\eta|\le R_*,
\qquad
|\zeta-\eta|\le R_*,
\]
and therefore $|\xi-\zeta|\le2R_*$. Furthermore,
\[
\mathcal J(\xi,\zeta)
\le\frac1{a_-}
\int_{B_{R_*}(\xi)}1\dd\eta
=\frac{\omega_dR_*^d}{a_-}.
\]
Consequently,
\[
j_-\ind_{\{|\xi-\zeta|<r_*/2\}}
\le\mathcal J(\xi,\zeta)
\le j_+\ind_{\{|\xi-\zeta|\le2R_*\}},
\]
where
\[
j_-:=\frac{k_*^2\omega_d(r_*/2)^d}{a_+}>0,
\qquad
j_+:=\frac{\omega_dR_*^d}{a_-}<\infty.
\]
All hypotheses of Lemma~\ref{lem:liouville} are satisfied.
Since $u$ is bounded and continuous, it follows that
$u\equiv c$ for some $c\in\R$.

The second identity in \eqref{eq:bipartite} then becomes
\[
0=b(\eta)v(\eta)
+c\int_{\R^d}\mathcal K(\xi,\eta)\dd\xi
=b(\eta)\bigl[v(\eta)+c\bigr].
\]
Since $b(\eta)>0$, we obtain $v=-c$. Thus
\[
A(\xi)=\frac{c_d}{2}+c,
\qquad
B(\eta)=\frac{c_d}{2}-c,
\]
and
\[
H(\xi,\eta)=c_d-\frac12|\xi-\eta|^2
=H_0(\xi,\eta).
\]
\end{proof}

\section{Uniform interior estimates in the original variables}
\label{sec:uniform}

The previous sections describe the unique limiting slack obtained after rescaling around
points of $K$. We now use this uniqueness to obtain estimates that hold uniformly over
all base points in $K$. We also express the estimates in the original
variables. The estimate for $r_\eps+s_\eps$ will be used in the next
section to compare two regularization parameters.

For $x\in K$, recall the definitions
\[
J_x:=\nabla T(x),
\qquad
h_\eps(x):=\left(\frac{\eps}{\lambda(x)}\right)^{1/(d+2)},
\]
from Section \ref{sec:rescaling} and set
\[
\eta_\eps(x,y):=
h_\eps(x)^{-1}J_x^{-1/2}[y-T(x)].
\]
Choosing the base point $p=x$ in \eqref{eq:coordinates}, we have $\xi=0$ and
$\eta=\eta_\eps(x,y)$. In particular,
\[
y=T(x)+h_\eps(x)J_x^{1/2}\eta_\eps(x,y).
\]

\begin{theorem}[Uniform interior slack and correction estimates]
\label{thm:uniform}
For every fixed $R>0$, we have
\begin{align}
\sup_{\substack{x\in K,\ y\in Y\\
|\eta_\eps(x,y)|\le R}}
\left|
\frac{H_\eps(x,y)}{h_\eps(x)^2}
-c_d+\frac12|\eta_\eps(x,y)|^2
\right|
&\longrightarrow0,
\label{eq:uniform-slack}\\
\sup_{\substack{x\in K,\ y\in Y\\
|\eta_\eps(x,y)|\le R}}
\left|
\frac{r_\eps(x)+s_\eps(y)}{h_\eps(x)^2}-c_d
\right|
&\longrightarrow0
\label{eq:uniform-correction}
\end{align}
as $\eps\downarrow0$.
\end{theorem}

\begin{proof}
Fix $R>0$. We start with \eqref{eq:uniform-slack}. Suppose it fails. Then there exist
$a>0$, a sequence $\eps_j\downarrow0$, and points
$p_j\in K$, $y_j\in Y$ such that, writing
\[
q_j:=T(p_j),
\qquad
J_j:= \nabla T(p_j),
\qquad
h_j:=h_{\eps_j}(p_j),
\qquad
\eta_j:=\eta_{\eps_j}(p_j,y_j),
\]
we have $|\eta_j|\le R$ and
\[
\left|
\frac{H_{\eps_j}(p_j,y_j)}{h_j^2}
-c_d+\frac12|\eta_j|^2
\right|\ge a
\qquad\text{for every }j.
\]
By the definition of $\eta_j$,
\[
y_j=q_j+h_jJ_j^{1/2}\eta_j.
\]
We now apply Proposition~\ref{prop:compactness} with base points $p_j$.
After passing to a subsequence, the rescaled corrections
\begin{align}\label{eq:A_j}
A_j(\xi)&:=
h_j^{-2}r_{\eps_j}(p_j+h_jJ_j^{-1/2}\xi),\\
B_j(\eta)&:=
h_j^{-2}s_{\eps_j}(q_j+h_jJ_j^{1/2}\eta)
\end{align}
converge locally uniformly to bounded Lipschitz functions $A,B$
satisfying \eqref{eq:frozen}. In particular,
\begin{align}\label{eq:help3}
|A_j(0)-A(0)|
+\sup_{|\eta|\le R}|B_j(\eta)-B(\eta)|
\longrightarrow0.
\end{align}
Proposition~\ref{prop:classification} gives
\begin{align}\label{eq:help1}
A(0)+B(\eta)=c_d
\qquad\text{for every }\eta\in\R^d.
\end{align}
Applying \eqref{eq:gap-taylor} with $\xi=0$ and $\eta=\eta_j$,
we obtain
\begin{align}\label{eq:help2}
\left|
h_j^{-2}D(p_j,y_j)-\frac12|\eta_j|^2
\right|
\le C_Rh_j.
\end{align}
Moreover, $h_j\to0$ by \eqref{eq:affine-data} and \eqref{eq:norm_bound}.
Using \eqref{eq:decomposition}, we therefore find
\begin{align*}
a
&\le
\left|
\frac{H_{\eps_j}(p_j,y_j)}{h_j^2}
-c_d+\frac12|\eta_j|^2
\right|\\
&\stackrel{\eqref{eq:decomposition}, \eqref{eq:A_j}}{=}
\left|
A_j(0)+B_j(\eta_j)-c_d
-\left[
h_j^{-2}D(p_j,y_j)-\frac12|\eta_j|^2
\right]
\right|\\
&\stackrel{\eqref{eq:help1},\eqref{eq:help2}}{\le}
|A_j(0)-A(0)|
+|B_j(\eta_j)-B(\eta_j)|
+C_Rh_j\\
&\hspace{0.52cm} \le
|A_j(0)-A(0)|
+\sup_{|\eta|\le R}|B_j(\eta)-B(\eta)|
+C_Rh_j
\stackrel{\eqref{eq:help3}}\to 0.
\end{align*}
This contradiction proves \eqref{eq:uniform-slack}.

We now prove \eqref{eq:uniform-correction}. For $x\in K$ and $y\in Y$ satisfying
$|\eta_\eps(x,y)|\le R$, Lemma~\ref{lem:rescaling}, applied with
base point $p=x$ and $\xi=0$ gives
\[
\left|
\frac{D(x,y)}{h_\eps(x)^2}
-\frac12|\eta_\eps(x,y)|^2
\right|
\le C_Rh_\eps(x).
\]
Since $r_\eps+s_\eps=H_\eps+D$ by \eqref{eq:decomposition},
\begin{align*}
\left|
\frac{r_\eps(x)+s_\eps(y)}{h_\eps(x)^2}-c_d
\right|
&\le
\left|
\frac{H_\eps(x,y)}{h_\eps(x)^2}
-c_d+\frac12|\eta_\eps(x,y)|^2
\right|\\
&\quad+
\left|
\frac{D(x,y)}{h_\eps(x)^2}
-\frac12|\eta_\eps(x,y)|^2
\right|.
\end{align*}
Taking the supremum, the first term tends
to zero by \eqref{eq:uniform-slack}, while the second is bounded by
\[
C_R\sup_{x\in K}h_\eps(x)
\le C_Rm^{-1/(d+2)}\eps^{1/(d+2)}
\to 0.
\]
This proves the claim.
\end{proof}

\section{Proof of Theorem \ref{thm:main}}
\label{sec:comparison}

The previous section gives an approximation for the slack that is
uniform over base points in $K$. We now use this approximation to
compare the slackes at two regularization parameters
$0<\delta\le\theta\eps$, where $0<\theta<1$ is fixed.

The comparison uses the decomposition
\[
H_\eps(x,y)=r_\eps(x)+s_\eps(y)-D(x,y).
\]
Since $D(x,y)$ is independent of the regularization parameter,
it cancels when we subtract $H_\delta(x,y)$ from $H_\eps(x,y)$.
The remaining terms are controlled by
Theorem~\ref{thm:uniform}. Applying Proposition~\ref{prop:localization} to
$(x,y)\in\spt\pi_\delta$ gives a bound on the rescaled
coordinates at both $\delta$ and $\eps$. This allows us to apply
Theorem~\ref{thm:uniform} to the same pair $(x,y)$ at both parameters. We first record the scaling behavior on the graph of $T$ and introduce a bound for the approximation errors. Taking $y=T(x)$ in \eqref{eq:uniform-slack} gives
\[
\sup_{x\in K}
\left|
\frac{H_\eps(x,T(x))}{h_\eps(x)^2}-c_d
\right|\longrightarrow0.
\]
By \eqref{eq:affine-data} and
$\gamma=2/(d+2)$ in \eqref{eq:constants},
\begin{align}\label{eq:hk}
h_\eps(x)^2=\lambda(x)^{-\gamma}\eps^\gamma.
\end{align}
Since $m\le\lambda(x)\le M$, we obtain
\begin{align*}
\sup_{x\in K}
\left|
\eps^{-\gamma}H_\eps(x,T(x))
-c_d\lambda(x)^{-\gamma}
\right|
\le
m^{-\gamma}
\sup_{x\in K}
\left|
\frac{H_\eps(x,T(x))}{h_\eps(x)^2}-c_d
\right|
\to 0.
\end{align*}
Thus
\[
H_\eps(x,T(x))
=c_d\lambda(x)^{-\gamma}\eps^\gamma
+o_K(\eps^\gamma),
\]
where the remainder divided by $\eps^\gamma$ tends uniformly
to zero for $x\in K$. In particular, for all sufficiently small
$\eps$,
\[
H_\eps(x,T(x))
\ge\frac{c_d}{2}M^{-\gamma}\eps^\gamma>0
\qquad\text{for every }x\in K.
\]
For the comparison at two parameters, define
\begin{equation}\label{eq:error}
e_\eps(x,y):=
\frac{r_\eps(x)+s_\eps(y)}{h_\eps(x)^2}-c_d.
\end{equation}
For fixed $R>0$ and sufficiently small $t>0$, set
\begin{equation}\label{eq:omega}
\omega_{K,R}(t):=
\sup_{\substack{0<\eps\le t,\ x\in K,\ y\in Y\\
|\eta_\eps(x,y)|\le R}}
|e_\eps(x,y)|.
\end{equation}
This function bounds the errors simultaneously for all
regularization parameters $\eps\le t$. Indeed, \eqref{eq:uniform-correction} implies that, for every
$a>0$, there exists $t_a>0$ such that
\[
0\le\omega_{K,R}(t)\le a
\qquad\text{whenever }0<t\le t_a.
\]
Consequently, $\omega_{K,R}$ is finite for sufficiently small $t$,
and $ \lim_{t\downarrow0}\omega_{K,R}(t)=0.$ Combining \eqref{eq:error} with \eqref{eq:decomposition}, we have
the exact identity
\begin{align}\label{eq:H_p}
H_\eps(x,y)
=c_dh_\eps(x)^2-D(x,y)+h_\eps(x)^2e_\eps(x,y),
\end{align}
with $|e_\eps(x,y)|\le\omega_{K,R}(t)$
whenever $0<\eps\le t$, $x\in K$, $y\in Y$, and
$|\eta_\eps(x,y)|\le R$.

\begin{proof}[Proof of Theorem~\ref{thm:main}]
Fix $0<\theta<1$. Let $0<\delta\le\theta\eps$ and $(x,y)\in\spt\pi_\delta, x\in K.$ We will choose an upper bound on $\eps$ below, independently of $\delta$ and $(x,y)$.

\emph{Step 1: bound the local coordinates at both scales.}
For $\eps$ sufficiently small, we also have $\delta\le\eps$
sufficiently small to apply Proposition~\ref{prop:localization}.
Thus $|y-T(x)|\le C_{\mathrm{loc}}\ell_\delta.$
Since $\|J_x^{-1/2}\|\le a_0^{-1/2}$ by \eqref{eq:hessian}, we obtain
\begin{align*}
|\eta_\delta(x,y)|
&=h_\delta(x)^{-1}
\bigl|J_x^{-1/2}[y-T(x)]\bigr|
\le
\frac{C_{\mathrm{loc}}\ell_\delta}
{\sqrt{a_0}\,h_\delta(x)}
\stackrel{\eqref{eq:comparison:hl}}{=}
\frac{C_{\mathrm{loc}}\lambda(x)^{1/(d+2)}}
{\sqrt{a_0}}
\le
\frac{C_{\mathrm{loc}}M^{1/(d+2)}}{\sqrt{a_0}}.
\end{align*}
Fix
\[
R>\frac{C_{\mathrm{loc}}M^{1/(d+2)}}{\sqrt{a_0}}.
\]
This choice only depends on the bounds in Assumption \ref{ass:data}.

At both parameters we use the same base point $x$ and hence the
same matrix $J_x$ and density factor $\lambda(x)$. Therefore
\[
\frac{h_\delta(x)}{h_\eps(x)}
=\left(\frac{\delta}{\eps}\right)^{1/(d+2)}
\]
and
\begin{align*}
\eta_\eps(x,y)
&=h_\eps(x)^{-1}J_x^{-1/2}[y-T(x)]
=\frac{h_\delta(x)}{h_\eps(x)}\eta_\delta(x,y).
\end{align*}
Since $\delta\le\eps$, it follows that $
|\eta_\eps(x,y)|
\le|\eta_\delta(x,y)|\le R.$
By \eqref{eq:omega}, both errors consequently satisfy
\begin{align}\label{eq:errorb_bound2}
|e_\eps(x,y)|\le\omega_{K,R}(\eps),
\qquad
|e_\delta(x,y)|\le\omega_{K,R}(\eps).
\end{align}

\emph{Step 2: subtract the slack identities.}
For the fixed pair $(x,y)$, write
\[
h:=h_\eps(x),
\qquad
k:=h_\delta(x),
\qquad
\omega:=\omega_{K,R}(\eps).
\]
By \eqref{eq:error} and \eqref{eq:decomposition},
\begin{align*}
H_\eps(x,y)
&=c_dh^2-D(x,y)+h^2e_\eps(x,y),\\
H_\delta(x,y)
&=c_dk^2-D(x,y)+k^2e_\delta(x,y).
\end{align*}
Subtracting these two equations gives
\begin{equation}\label{eq:subtract}
H_\eps(x,y)
=H_\delta(x,y)+c_d(h^2-k^2)
+h^2e_\eps(x,y)-k^2e_\delta(x,y).
\end{equation}
By Lemma~\ref{lem:support},
\[
(x,y)\in\spt\pi_\delta
\quad\Longrightarrow\quad
H_\delta(x,y)\ge0.
\]
Using \eqref{eq:errorb_bound2}, we therefore obtain
\[
H_\eps(x,y)
\ge c_d(h^2-k^2)-\omega(h^2+k^2).
\]

Set
\[
\sigma:=\frac{k^2}{h^2}
\stackrel{\eqref{eq:hk}}{=}\left(\frac{\delta}{\eps}\right)^\gamma.
\]
Since $\delta\le\theta\eps$, we have
$0<\sigma\le\theta^\gamma<1$. Consequently,
\begin{align}\label{eq:step2}
H_\eps(x,y)
&\ge h^2\bigl[c_d(1-\sigma)-\omega(1+\sigma)\bigr]\nonumber\\
&\ge h^2\bigl[c_d(1-\theta^\gamma)
-\omega(1+\theta^\gamma)\bigr].
\end{align}

\emph{Step 3: choosing a threshold for $\eps$.}
Since $\theta<1$,
\[
\frac{c_d(1-\theta^\gamma)}
{2(1+\theta^\gamma)}>0.
\]
Moreover, $\omega_{K,R}(\eps)\to0$ as $\eps\downarrow0$.
We may therefore choose $\eps_{K,\theta}>0$ such that
Proposition \ref{prop:localization} applies at every parameter below $\eps_{K,\theta}$
and
\[
\omega_{K,R}(\eps)
\le
\frac{c_d(1-\theta^\gamma)}
{2(1+\theta^\gamma)}
\qquad\text{for all }0<\eps<\eps_{K,\theta}.
\]
Equation \eqref{eq:step2} then gives
\[
H_\eps(x,y)
\ge\frac{c_d}{2}(1-\theta^\gamma)h_\eps(x)^2.
\]
Finally,
\[
h_\eps(x)^2
\stackrel{\eqref{eq:affine-data}}{=}\lambda(x)^{-\gamma}\eps^\gamma
\ge M^{-\gamma}\eps^\gamma,
\]
and hence
\[
H_\eps(x,y)
\ge
\frac{c_d}{2}M^{-\gamma}
(1-\theta^\gamma)\eps^\gamma>0.
\]
This proves \eqref{eq:main}. Applying
Lemma~\ref{lem:support} once more, we conclude that
\[
\spt\pi_\delta\cap(K\times Y)
\subseteq
\spt\pi_\eps\cap(K\times Y)
\]
whenever
\[
0<\eps<\eps_{K,\theta},
\qquad
0<\delta\le\theta\eps.
\]
\end{proof}

The threshold $\eps_{K,\theta}$ may depend on $K$, $\theta$ and the bounds in Assumption \ref{ass:data}. The coefficient
\[
\frac{c_d}{2}M^{-\gamma}(1-\theta^\gamma)
\]
in the slack lower bound is independent of $K$.

\section{Acknowledgements}

JW acknowledges support by Novo Nordisk and Villum Fonden.
The authors used ChatGPT (OpenAI) to assist with proof reading, image generation and references. They take full responsibility for the final manuscript.

\appendix
\section{Proof of Lemma \ref{lem:liouville}}
\label{app:liouville}

We now prove Lemma~\ref{lem:liouville}, which says the following: if
\[
\int_{\R^d}\mathcal J(x,y)[u(x)-u(y)]\dd y=0,
\]
where $\mathcal J$ is symmetric and satisfies
\eqref{eq:Jbounds} and $u$ is bounded and continuous, then $u$ is constant.

The proof first expresses the equation as an averaging identity.
We then verify the hypotheses of the elliptic Harnack theorem
\cite[Theorem~6.3]{MSC2015}. This theorem compares the supremum
and infimum of a nonnegative harmonic function on a ball, with
a constant independent of sufficiently large radii. Applying
this estimate to $u-\inf_{\R^d}u$ and letting the radius tend
to infinity will prove the result.

\begin{proof}[Proof of Lemma~\ref{lem:liouville}]
We proceed in four steps. In the first three steps, we verify the hypotheses of
\cite[Theorem~6.3]{MSC2015}. In the last step, we apply the Harnack inequality to conclude.

\emph{Step 1: express \eqref{eq:u} as an average.}
Let us define
\[
w(x):=\int_{\R^d}\mathcal J(x,y)\dd y.
\]
The bounds \eqref{eq:Jbounds} then imply
\begin{equation}\label{eq:liouville-weights}
0<w_-:=j_-\omega_dr^d
\le w(x)
\le j_+\omega_dR_0^d=:w_+<\infty.
\end{equation}
Furthermore, \eqref{eq:u} is equivalent to
\[
u(x)=Pu(x),
\qquad
Pf(x):=\frac1{w(x)}
\int_{\R^d}\mathcal J(x,y)f(y)\dd y.
\]
As $P1=1$ by definition of $w$, we think of $Pf$ as an average of $f$. Indeed, defining the measure $\dd\vartheta(x):=w(x)\dd x$
and the density $p(x,y):=\frac{\mathcal J(x,y)}{w(x)w(y)}$
we have
\[
Pf(x)=\int_{\R^d}p(x,y)f(y)\dd\vartheta(y),
\]
and
\[
p(x,y)=p(y,x),
\qquad
\int_{\R^d}p(x,y)\dd\vartheta(y)=1.
\]
Thus $P$ is reversible with respect to $\vartheta$.
The density $p$ is the transition density relative to
$\vartheta$ used below. In particular, by \eqref{eq:Jbounds} and \eqref{eq:liouville-weights},
\begin{align}\label{eq:liouville-pbounds}
a\ind_{\{|x-y|<r\}}
\le p(x,y)
\le b\ind_{\{|x-y|\le R_0\}},
\end{align}
for $a:=j_-/w_+^2$ and $b:=j_+/w_-^2.$

\emph{Step 2: verify volume doubling and the Poincar\'e inequality.}
Write
\[
V(x,t):=\vartheta(B_t(x)).
\]
By \eqref{eq:liouville-weights},
\begin{equation}\label{eq:liouville-volumes}
w_-\omega_dt^d
\le V(x,t)
\le w_+\omega_dt^d.
\end{equation}
Consequently,
\[
V(x,2t)
\le 2^d\frac{w_+}{w_-}V(x,t)
\qquad\text{for every }x\in\R^d,\ t>0.
\]
This is the volume doubling property, with a constant independent
of $t$ and $x$.

We next prove the Poincar\'e inequality required for the Harnack
theorem. Fix $h:= r/2$.
For $f\in L^2_{\mathrm{loc}}(\R^d)$, define
\[
G_hf(x):=
\frac1{V(x,h)}
\int_{B_h(x)}|f(x)-f(y)|^2\dd\vartheta(y).
\]
The quantity $G_hf(x)$ measures the average squared difference between
$f(x)$ and values of $f$ within distance $h$. 

Next, we choose a standard mollifier
\[
\zeta\in C_c^\infty(B_{1/2}(0)),
\qquad
\zeta\ge0,
\qquad
\int_{\R^d}\zeta=1,
\]
and set
\[
\zeta_h(z):=h^{-d}\zeta(z/h),
\qquad
f_h:=f*\zeta_h.
\]
Jensen's inequality gives, for almost every $x$,
\begin{align*}
|f(x)-f_h(x)|^2
&=
\left|
\int_{\R^d}\zeta_h(x-y)[f(x)-f(y)]\dd y
\right|^2\\
&\le
\int_{\R^d}\zeta_h(x-y)|f(x)-f(y)|^2\dd y\\
&\le
Ch^{-d}\int_{B_h(x)}|f(x)-f(y)|^2\dd y.
\end{align*}
Moreover, $\int_{\R^d}\nabla\zeta_h=0$, so
\[
\nabla f_h(x)
=
\int_{\R^d}\nabla\zeta_h(x-y)[f(y)-f(x)]\dd y.
\]
Since
\[
\int_{\R^d}|\nabla\zeta_h|^2
=h^{-d-2}\int_{\R^d}|\nabla\zeta|^2,
\]
Cauchy--Schwarz yields
\[
h^2|\nabla f_h(x)|^2
\le
Ch^{-d}\int_{B_h(x)}|f(x)-f(y)|^2\dd y.
\]
Combining these estimates, we obtain
\begin{equation}\label{eq:liouville-smoothing}
|f(x)-f_h(x)|^2+h^2|\nabla f_h(x)|^2
\le
Ch^{-d}\int_{B_h(x)}|f(x)-f(y)|^2\dd y.
\end{equation}

Now we fix $x_0\in\R^d$ and $t\ge h$, and write $B_t:=B_t(x_0)$.
Define the two means
\[
f_{B_t,\vartheta}
:=
\frac1{\vartheta(B_t)}\int_{B_t}f\dd\vartheta,
\qquad
c_t:=\frac1{|B_t|}\int_{B_t}f_h\dd x.
\]
The $\vartheta$-mean minimizes the squared distance to constants.
Using this fact, \eqref{eq:liouville-weights}, and the ordinary
Euclidean Poincar\'e inequality for $f_h$, we obtain
\begin{align*}
\int_{B_t}|f-f_{B_t,\vartheta}|^2\dd\vartheta
&\le
\int_{B_t}|f-c_t|^2\dd\vartheta\\
&\le
2w_+\int_{B_t}|f-f_h|^2\dd x
+2w_+\int_{B_t}|f_h-c_t|^2\dd x\\
&\le
2w_+\int_{B_t}|f-f_h|^2\dd x
+Cw_+t^2\int_{B_t}|\nabla f_h|^2\dd x\\
&\stackrel{\eqref{eq:liouville-smoothing}}{\le}
Cw_+\left(\frac th\right)^2h^{-d}
\int_{B_t}\int_{B_h(x)}
|f(x)-f(y)|^2\dd y\dd x,
\end{align*}
where we have enlarged the constant $C$ such that both terms have the common factor $(t/h)^2\ge1$.
To express the right-hand side in terms of $\vartheta$, observe
from \eqref{eq:liouville-weights} and
\eqref{eq:liouville-volumes} that
\[
G_hf(x)
\ge
\frac{w_-}{w_+\omega_dh^d}
\int_{B_h(x)}|f(x)-f(y)|^2\dd y.
\]
Taking the outer integral wrt.~$\vartheta$, it  follows that
\begin{equation}\label{eq:liouville-poincare}
\int_{B_t}|f-f_{B_t,\vartheta}|^2\dd\vartheta
\le
C\left(\frac{w_+}{w_-}\right)^2
\left(\frac th\right)^2
\int_{B_t}G_hf(x)\dd\vartheta(x).
\end{equation}
Since $h$ is fixed, this has the required form
\[
\int_{B_t}|f-f_{B_t,\vartheta}|^2\dd\vartheta
\le C_ht^2\int_{B_t}G_hf\dd\vartheta,
\]
where $C_h$ is independent of $x_0$, $t$, and $f$.

Finally, straight segments join any two points in $\R^d$.
Subdividing such a segment gives steps of length at most $h/2$
whose total length equals the distance between the endpoints.
Thus the geometric hypothesis of the Harnack theorem holds
with $C=1$ and a scale $b$ strictly smaller than $h$.

\emph{Step 3: construct an averaging operator satisfying the
remaining kernel condition.}
The Harnack theorem also requires the transition density $q$
of the averaging operator $\mathcal{P}$ (defined below) to satisfy
\[
q_2(x,y)\ge\alpha q(x,y)
\]
for some $\alpha>0$, where
\[
q_2(x,y):=
\int_{\R^d}q(x,z)q(z,y)\dd\vartheta(z)
\]
is the density corresponding to two successive applications
of this operator. We now construct the averaging operator $\mathcal{P}$ from $P$. Let $p_i$ denote the
density of $P^i$ relative to $\vartheta$, i.e.~
\begin{align}\label{eq:pee}
p_1=p,
\qquad
p_{i+1}(x,y)
=
\int_{\R^d}p_i(x,z)p(z,y)\dd\vartheta(z).
\end{align}
By induction, repeated integration and symmetry yield
\begin{align}\label{eq:prop_pee}
p_i(x,y)=p_i(y,x),
\qquad
\int_{\R^d}p_i(x,y)\dd\vartheta(y)=1.
\end{align}
Moreover,
\[
p_{i+1}(x,y)
\stackrel{\eqref{eq:liouville-pbounds}}{\le} b\int_{\R^d}p_i(x,z)\dd\vartheta(z)=b,
\]
and, again by \eqref{eq:liouville-pbounds},
\[
p_i(x,y)=0
\qquad\text{if }|x-y|>iR_0.
\]
Now choose an integer
\[
N\ge\max\left\{2,\left\lceil\frac{4R_0}{r}\right\rceil\right\}.
\]
For $|x-y|\le R_0$, set
\[
x_i:=x+\frac{i}{N+1}(y-x),
\qquad i=0,\ldots,N+1.
\]
We now bound $p_{N+1}(x,y)$ from below. For this, in the integral in \eqref{eq:pee} defining $p_{N+1}(x,y)$, restrict the intermediate
points to
\[
z_i\in B_{r/8}(x_i),
\qquad i=1,\ldots,N,
\]
and put $z_0:=x$, $z_{N+1}:=y$. Then
\[
|z_{i+1}-z_i|
\le\frac{R_0}{N+1}+\frac r4<r
\qquad i=0,\ldots,N.
\]
By \eqref{eq:liouville-pbounds}, every factor
$p(z_i,z_{i+1})$ is at least $a$. Consequently,
\begin{align*}
p_{N+1}(x,y)
&\ge
a^{N+1}\prod_{i=1}^N
\vartheta(B_{r/8}(x_i))\\
&\stackrel{\eqref{eq:liouville-weights}}{\ge}
a^{N+1}
\bigl[w_-\omega_d(r/8)^d\bigr]^N
=:c_*>0.
\end{align*}
Since $p\le b$ and $p(x,y)=0$ for $|x-y|>R_0$ by \eqref{eq:liouville-pbounds}, we conclude that
\begin{equation}\label{eq:liouville-power-lower}
p_{N+1}(x,y)\ge\beta p(x,y),
\qquad
1\ge \beta:=\min\{1,c_*/b\}>0.
\end{equation}
We are now in a position to define the averaging operator and its transition density
\[
\mathcal P:=\frac1N\sum_{i=1}^NP^i,
\qquad
q(x,y):=\frac1N\sum_{i=1}^Np_i(x,y).
\]
Recalling \eqref{eq:prop_pee}, $q$ is symmetric, has integral one in its second variable,
and satisfies
\begin{align}\label{eq:qbounds}
\frac aN\ind_{\{|x-y|<r\}}
\le q(x,y)
\le b\ind_{\{|x-y|\le NR_0\}}.
\end{align}
Its two-step density satisfies
\begin{align}\label{eq:liouville-q2}
\begin{split}
q_2(x,y)
&=\frac1{N^2}\sum_{i,j=1}^Np_{i+j}(x,y)\\
&\ge
\frac1{N^2}\sum_{k=2}^{N+1}p_k(x,y)\\
&= \frac1{N^2} \Big(p_{N+1} +\sum_{k=2}^{N}p_k(x,y) \Big)\\
&\stackrel{\eqref{eq:liouville-power-lower}}{\ge} \frac1{N^2} \Big( \beta p_1 +\beta \sum_{k=2}^{N}p_k(x,y) \Big)\\
&=
\frac\beta{N^2}\sum_{k=1}^Np_k(x,y)\\
&=\frac\beta Nq(x,y).
\end{split}
\end{align}
For completeness, the volume-normalized kernel bounds required
by the Harnack theorem follow from rewriting \eqref{eq:qbounds} via
\eqref{eq:liouville-volumes}:
\[
\frac{c_1}{V(x,h)}\ind_{\{|x-y|\le h\}}
\le q(x,y)
\le
\frac{C_1}{V(x,NR_0)}\ind_{\{|x-y|\le NR_0\}},
\]
where one may take
\[
c_1:=\frac aNw_-\omega_dh^d,
\qquad
C_1:=bw_+\omega_d(NR_0)^d.
\]
Together with symmetry, normalization, and
\eqref{eq:liouville-q2}, these are precisely the required
compatibility conditions.

Finally, $Pu=u$ implies $P^iu=u$ for every $i\ge1$, and hence
\[
\mathcal Pu
=\frac1N\sum_{i=1}^NP^iu=u.
\]

\emph{Step 4: apply Harnack on expanding balls.}
The preceding steps verify the hypotheses of
\cite[Theorem~6.3]{MSC2015}. Thus there exist
$c\in(0,1)$, $C_H<\infty$, and $t_0>0$ such that every bounded
nonnegative $\mathcal P$-harmonic function $f$ satisfies
\begin{equation}\label{eq:liouville-harnack}
\mathop{\rm ess\,sup}_{B_{ct}(0)}f
\le
C_H\mathop{\rm ess\,inf}_{B_{ct}(0)}f
\qquad\text{for every }t\ge t_0.
\end{equation}
The constant $C_H$ is independent of $t$. We now set
\[
m:=\inf_{\R^d}u,
\qquad
\widetilde u:=u-m.
\]
Since $u$ is bounded, $m$ is finite. Moreover,
\[
\widetilde u\ge0,
\qquad
\inf_{\R^d}\widetilde u=0,
\qquad
\mathcal P\widetilde u
=\mathcal Pu-m\mathcal P1=u-m=\widetilde u.
\]
The function $\widetilde u$ is continuous, and $\vartheta$ is
comparable to Lebesgue measure by
\eqref{eq:liouville-weights}. Its essential and ordinary extrema
on balls therefore agree. Applying
\eqref{eq:liouville-harnack} gives
\[
\sup_{B_{ct}(0)}\widetilde u
\le C_H\inf_{B_{ct}(0)}\widetilde u
\qquad\text{for every }t\ge t_0.
\]
As $t\to\infty$, these balls exhaust $\R^d$, so
\[
\inf_{B_{ct}(0)}\widetilde u
\to\inf_{\R^d}\widetilde u=0,
\qquad
\sup_{B_{ct}(0)}\widetilde u
\longrightarrow\sup_{\R^d}\widetilde u.
\]
Consequently,
\[
0\le\sup_{\R^d}\widetilde u
\le C_H\inf_{\R^d}\widetilde u=0.
\]
Thus $\widetilde u=0$, and $u= m$.
\end{proof}

\bibliographystyle{plain}
\bibliography{refs}

\begin{thebibliography}{10}

\bibitem{BlondelSeguyRolet2018}
Mathieu Blondel, Vivien Seguy, and Antoine Rolet.
\newblock Smooth and sparse optimal transport.
\newblock In {\em Proceedings of the Twenty-First International Conference on
  Artificial Intelligence and Statistics}, volume~84 of {\em Proceedings of
  Machine Learning Research}, pages 880--889. PMLR, 2018.

\bibitem{ChengKochXiao2026GlobalRegularity}
Xiaopeng Cheng, Lukas Koch, and Haotian Xiao.
\newblock Global regularity and sparsity for regularised optimal transport,
  2026.
\newblock arXiv:2610.01928.

\bibitem{EcksteinNutz2024Quantization}
Stephan Eckstein and Marcel Nutz.
\newblock Convergence rates for regularized optimal transport via quantization.
\newblock {\em Mathematics of Operations Research}, 49(2):1223--1240, 2024.

\bibitem{EssidSolomon2018}
Montacer Essid and Justin Solomon.
\newblock Quadratically regularized optimal transport on graphs.
\newblock {\em SIAM Journal on Scientific Computing}, 40(4):A1961--A1986, 2018.

\bibitem{GarrizMolinaGonzalezSanzMordant2024Infinitesimal}
Alejandro Garriz-Molina, Alberto Gonz{\'a}lez-Sanz, and Gilles Mordant.
\newblock Infinitesimal behavior of quadratically regularized optimal transport
  and its relation with the {Porous Medium Equation}, 2024.
\newblock arXiv:2407.21528.

\bibitem{GonzalezSanzGvalaniKoch2026LocalSparsity}
Alberto Gonz{\'a}lez-Sanz, Rishabh~S. Gvalani, and Lukas Koch.
\newblock Sharp local sparsity of regularized optimal transport, 2026.
\newblock arXiv:2604.00843.

\bibitem{GonzalezSanzNutz2025LinearPrograms}
Alberto Gonz{\'a}lez-Sanz and Marcel Nutz.
\newblock Quantitative convergence of quadratically regularized linear
  programs.
\newblock {\em Applied Mathematics \& Optimization}, 91(3), 2025.
\newblock Article 68.

\bibitem{GonzalezSanzNutz2026GeometryI}
Alberto Gonz{\'a}lez-Sanz and Marcel Nutz.
\newblock Geometry and convergence of quadratically regularized optimal
  transport {I}, 2026.
\newblock arXiv:2609.20400.

\bibitem{GonzalezSanzNutz2026GeometryII}
Alberto Gonz{\'a}lez-Sanz and Marcel Nutz.
\newblock Geometry and convergence of quadratically regularized optimal
  transport {II}, 2026.
\newblock arXiv:2610.02077.

\bibitem{gonzalez2026sparsity}
Alberto Gonz{\'a}lez-Sanz and Marcel Nutz.
\newblock Sparsity of quadratically regularized optimal transport: Scalar case.
\newblock {\em SIAM Journal on Mathematical Analysis}, 58(4):4162--4183, 2026.

\bibitem{gonzalez2024monotonicity}
Alberto Gonz{\'a}lez-Sanz, Marcel Nutz, and Andr{\'e}s~Riveros Valdevenito.
\newblock Monotonicity in quadratically regularized linear programs.
\newblock {\em arXiv preprint arXiv:2408.07871}, 2024.

\bibitem{GvalaniKoch2026UniformRegularity}
Rishabh~S. Gvalani and Lukas Koch.
\newblock Sparsity and uniform regularity for regularised optimal transport,
  2026.
\newblock arXiv:2601.05130.

\bibitem{LorenzMahler2022}
Dirk Lorenz and Hinrich Mahler.
\newblock {Orlicz} space regularization of continuous optimal transport
  problems.
\newblock {\em Applied Mathematics \& Optimization}, 85, 2022.
\newblock Article 14.

\bibitem{LorenzMannsMeyer2021}
Dirk~A. Lorenz, Paul Manns, and Christian Meyer.
\newblock Quadratically regularized optimal transport.
\newblock {\em Applied Mathematics \& Optimization}, 83(3):1919--1949, 2021.

\bibitem{MSC2015}
Mathav Murugan and Laurent Saloff-Coste.
\newblock Harnack inequalities and {Gaussian} estimates for random walks on
  metric measure spaces, 2015.
\newblock Preprint, arXiv:1506.07539v1,
  \url{https://arxiv.org/abs/1506.07539v1}.

\bibitem{Nutz2025}
Marcel Nutz.
\newblock Quadratically regularized optimal transport: Existence and
  multiplicity of potentials.
\newblock {\em SIAM Journal on Mathematical Analysis}, 57(3):2622--2649, 2025.
\newblock Preprint version: \url{https://arxiv.org/abs/2404.06847v3}.

\bibitem{wiesel2025sparsity}
Johannes Wiesel and Xingyu Xu.
\newblock Sparsity of quadratically regularized optimal transport: Bounds on
  concentration and bias.
\newblock {\em SIAM Journal on Mathematical Analysis}, 57(6):6498--6521, 2025.

\end{thebibliography}

\end{document}